\documentclass[12pt]{amsart}
\usepackage{amsmath,amsfonts,amstext,amsthm,amssymb,pxfonts}
\usepackage[colorlinks,citecolor=red,hypertexnames=false]{hyperref}

\theoremstyle{plain}

\newtheorem{theorem}[equation]{Theorem} 
 
\newtheorem{conjecture}[equation]{Conjecture}

\theoremstyle{definition}

\theoremstyle{remark}
\newtheorem{remark}[equation]{Remark}

\allowdisplaybreaks

\numberwithin{equation}{section}

\title{A Remark On coefficients 
of Jacobi matrices arising from a Schr{\"o}dinger operator}
\author[ A.~Vagharshakyan ] {Armen Vagharshakyan }
\address{\newline Armen Vagharshakyan \newline Mathematics Department, 
\newline Brown University, \newline 151 Thayer St, \newline Providence, RI 02912 USA}
\email{armen@math.brown.edu}
\begin{document}

\begin{abstract}
A discrete analogue of a Schr{\"o}dinger type operator proposed by J. Bellissard \cite{BBM} has a singular continuous spectrum. 
 In this remark we answer the conjecture formulated in \cite{BMM} on the coefficients of that operator. It turns out that the coefficients have a more complicated behavior than it was conjectured.
\end{abstract}
\maketitle
\noindent
\newline
\textit{AMS subject classification:} 47B36 (Jacobi (tridiagonal) operators and generalizations), 47A25 (Spectral sets)
\newline
\textit{PACS Numbers:}  05.50.+q (Lattice theory and statistics (Ising, Potts, etc.))
\section{\textbf{Introduction}}
\subsection{Dyson's example}
Consider the following discrete analog of the 
Schr{\"o}dinger operator:
\begin{equation*}
(Hu)(n)=p(n)u(n+1)+q(n)u(n-1)+V(n)u(n).
\end{equation*}
acting in the spaces $l_2(Z)$ or
$l_2(N)$.

	For illustrative purposes we describe the example of F. Dyson (see \cite{D}) 
where such an operator naturally arises.

Consider a chain of $N$ masses, each coupled to 
its nearest neighbors by elastic springs obeying Hooke's law.
Let all motions take place in one dimension 
so that each mass is described by a single coordinate.

Let particle number $j$ in the chain have mass $m_j$, and 
let its displacement from its equilibrium position be $x_j$. Let the elastic 
modules of the spring between particles $j$ and $j+1$ 
be $K_j$. Then the equations of the system's motion are:
\begin{equation*}
m_j\frac{d^2x_j}{dt^2}=K_j(x_{j+1}-x_j)+K_{j-1}(x_{j-1}-x_j),\quad j=2,3,\dots,N-1.
\end{equation*}

On the particles $x_1,\,\, x_N$ situated on the ends
we put certain boundary conditions.

It is convenient to introduce new variables:
\begin{equation*}
y_j=x_j\sqrt{m_j}
\end{equation*}
and new constants
$\lbrace
\lambda_j
\rbrace_{j=1}^{2N-2} $
given by:
\begin{equation*}
\lambda_{2j-1}=\frac{K_j}{m_j},\quad 
\lambda_{2j}=\frac{K_j}{m_{j+1}}.
\end{equation*}
Then the equations of motion take the symmetric form:
\begin{equation*}
\frac{d^2y_j}{dt^2}=y_{j+1}\sqrt{\lambda_{2j-1}\lambda_{2j}}+
y_{j-1}\sqrt{\lambda_{2j-3}\lambda_{2j-2}}-
(\lambda_{2j-1}+\lambda_{2j-2})y_j,\quad j=2,3,\dots,N-1.
\end{equation*}

For each mode 
\begin{equation*}
y_j(t)=u(j)e^{i E t}
\end{equation*}
we have
\begin{equation*}
u(j+1)\sqrt{\lambda_{2j-1}\lambda_{2j}}+
u(j-1)\sqrt{\lambda_{2j-3}\lambda_{2j-2}}-
(\lambda_{2j-1}+\lambda_{2j-2})u(j)=-E^2u(j),
\end{equation*}
where $j=2,3,\dots,N-1$.

For physical applications it is important to understand 
the properties of the spectrum, i.e. the properties of values of $E$ for 
which the equations have a nontrivial solution.

	In relation with this problem, mathematically it is more convenient 
to investigate the infinite dimensional case.

\subsection{Schr{\"o}dinger operators}

	Consider the following two types of operators:
\begin{equation}\label{T1}
T_1(u)(j)=u(j+1)\sqrt{\lambda_{2j-1}\lambda_{2j}}+
u(j-1)\sqrt{\lambda_{2j-3}\lambda_{2j-2}}-
(\lambda_{2j-1}+\lambda_{2j-2})u(j),\quad j\in Z
\end{equation}
where $u\in l_2(Z)$;
and
\begin{equation}\label{T2}
T_2(u)(j)=u(j+1)\sqrt{\lambda_{2j-1}\lambda_{2j}}+
u(j-1)\sqrt{\lambda_{2j-3}\lambda_{2j-2}}-
(\lambda_{2j-1}+\lambda_{2j-2})u(j),\quad j\in N
\end{equation}
where $u\in l_2(N)$.

	Note that for the first type \eqref{T1} there is no boundary condition while 
for the second type \eqref{T2} we take into account the existence of the boundary, too.
  
A typical equation of the first type is the almost Mathieu equation:
\begin{equation}\label{almostmathieu}
u(j+1)+u(j-1)+\lambda u(j)\cos2\pi(\theta-\alpha j)=Eu(j),\quad j\in Z
\end{equation}
where $u\in l_2(Z)$ and $\alpha$ is an irrational number. 

	A. Avron and S. Jitomirskaya (see \cite{AJ}) proved:
\begin{theorem}
Let $\alpha$ be irrational number and $\lambda\neq 0$. Then for almost all $\theta$ the spectrum
of the operator
\begin{equation*}
T(u)(j)=u(j+1)+u(j-1)+\lambda\cos2\pi\left(\theta-j\alpha\right)u(j),\quad j\in Z
\end{equation*}
is singular continuous.
\end{theorem}

	A typical equation of the second type is:

\begin{equation}\label{secondtypetypical}
(Hu)(j)=u(j+1)+u(j-1)+V(j) u(j),\quad j\in N
\end{equation}
where $u\in l_2(N)$. 

	Following \cite[p. 212]{CFKS} let's point out that general wisdom used 
to say that Schr{\"o}dinger operators should have absolutely 
continuous spectrum plus some discrete point spectrum, while
singular continuous spectrum is a pathology that should not occur in examples with $V$ bounded.
This general picture was proven to be wrong by Pearson (see \cite{P1}, \cite{P2})
who constructed a potential $V$ such that the operator
$H=H_0+V$ has singular continuous spectrum. His potential
$V$ consists of bumps further and further apart with the height of bumps possibly decreasing. 
Furthermore, singular continuous spectrum occurs in the innocent - looking almost
Mathieu equation \eqref{almostmathieu}.

\subsection{Bellissard's example}

One would like to find a potential function whose properties resemble physical phenomena closer.
For example, it would be more preferable to have an almost periodic potential.

	With this aim let us consider an operator of the form  \eqref{T2} proposed in
\cite{BBM}:
\begin{equation*}
(Hu)(j)=\sqrt{R_{j+1}}u(j+1)+\sqrt{R_j}u(j-1),\quad j\in N
\end{equation*}
where the numbers $R_n$ are defined recursively by \eqref{r1},\eqref{r2},\eqref{r3}.

	In paper \cite[p. 134]{BMM} the following results are proved for $\lambda>2$:
	
Proposition 1.\begin{equation*}
0<R_{2n}<R_n,\quad   0<R_{2n}\leq 1,
\end{equation*}

Proposition 2.  
\begin{equation*}
\lim_{k\to \infty}R_{p2^k}=0,
\end{equation*}

Proposition 3.  
\begin{equation*}
\lim_{k\to \infty}R_{p2^k+s}=R_s.
\end{equation*}

	Let us note that if in Proposition 3 the limit is uniform for $p$ and $s$ then
$R_n$ is a limit - periodic sequence (see \cite{Bohr}). 

	In ( \cite{BMM}, p. 135) it was conjectured that the sequence $R_n$ splits in the following way:
$R_{k2^r+n}$ lies between $R_n$ and $R_{2^r+n}$. The authors further point out that would their conjecture be true then the set $\lbrace R_n\rbrace$ would be perfect.

In this paper we prove that for $\lambda>2$ the sequence $R_n$ splits 
into 4 parts (see \eqref{t.main}). However, there is no 
further splitting (see \eqref{r.furthersplitting}). 

This shows that the structure of the set $R_n$ is much more complicated than it was conjectured in \cite{BMM}. 

\section{\textbf{Formulation of the Results}}
Let $\lambda>2$. Corresponding to 
$\lambda$ let us  discuss the  numerical sequence $R_n$ (where $n=0,1,\dots$) 
defined recursively by:
\begin{equation}\label{r1}
R_0=0
\end{equation}
\begin{equation}\label{r2}
R_{2n}+R_{2n+1}=\lambda
\end{equation}
\begin{equation}\label{r3}
R_{2n} R_{2n-1}=R_n
\end{equation}
According to \cite[p.~135]{BMM} numerical studies show that the set 
$\lbrace R_n \rbrace_{n=0}^{\infty}$ separates into disjoint subsets as follows:
\begin{conjecture}\label{c.main}
Let $\lambda>2$ and $R_n$ be defined by \eqref{r1},\eqref{r2},\eqref{r3} then for 
any  $n=0,1,2,\dots$ the following inequalities hold:
\begin{equation}\label{c1}
R_0\leq  R_{4n} \leq R_4,
\end{equation} 
\begin{equation}\label{c2}
R_6\leq R_{4n+2}\leq R_2,
\end{equation}
\begin{equation}\label{c3}
R_3\leq R_{4n+3} \leq R_7,
\end{equation} 
and
\begin{equation}\label{c4}
R_5\leq R_{4n+1}\leq R_1
\end{equation}
\end{conjecture}
In the light of this conjecture we prove the following theorem:
\begin{theorem}\label{t.main}
Let $\lambda>2$ and $R_n$ be defined by \eqref{r1},\eqref{r2},\eqref{r3} then for any 
$n=0,1,2,\dots$ the following inequalities hold:
\begin{equation}\label{p1}
0<R_{4n}\leq \frac{1}{\lambda -1}
\end{equation}
\begin{equation}\label{p2}
1-\frac{1}{\lambda -1}\leq R_{4n+2}<1
\end{equation}
\begin{equation}\label{p3}
\lambda-1<R_{4n+3}\leq \lambda-1+\frac{1}{\lambda-1}
\end{equation}
\begin{equation}\label{p4}
\lambda-\frac{1}{\lambda-1}\leq R_{4n+1}<\lambda
\end{equation}
\end{theorem}

\begin{remark}
Theorem \eqref{t.main} compares with the conjecture \eqref{c.main} in the following way:
\\
1. the bounds in the inequalities \eqref{c1} and \eqref{p1} are the same and obviously are sharp; 
\\
2. the bounds in the inequalities \eqref{c4} and \eqref{p4} are the same and obviously are sharp;
\\
3. the inequality \eqref{p3} proved in \eqref{t.main} is sharper than the inequality \eqref{c3} 
conjectured in \eqref{c.main}, indeed one can calculate that for $\lambda>2$ we have:
\begin{equation*}
R_7=\frac{\lambda^3-2\lambda^2+\lambda-1}{\lambda^2-\lambda-1}>\lambda-1+\frac{1}{\lambda-1}
\end{equation*}
\\
4. the inequality \eqref{p2} proved in theorem \eqref{t.main} is weaker than the 
inequality \eqref{c2} conjectured in \eqref{c.main}, indeed one can calculate that for $\lambda>2$ we have:
\begin{equation*}
R_6=\frac{(\lambda-1)^2}{\lambda^2-\lambda-1}>1-\frac{1}{\lambda-1}
\end{equation*}
But as it turns out the conjectured inequality \eqref{c2} is not true, indeed for $\lambda=2.1$ we have $R_{10}<R_6$.
Interestingly, for large values of $\lambda$ the lower bound of \eqref{c2} seems to be true.
\end{remark}
\begin{remark}\label{r.furthersplitting}
In \cite[p.~135]{BMM} it is further conjectured that $R_{k2^r+n}$ lies
 between $R_n$ and $R_{2^r+n}$. The same example with $\lambda=2.1$ and 
 $R_{10}<R_6<R_2$ comes to prove that this conjecture is not true.
\end{remark}
\section{Proof of Theorem \eqref{t.main}}
\subsection{Step 1}
From the recurrent formula \eqref{r2} it follows that:
\begin{equation}\label{additiveimplies}
\begin{split}
R_{4n}+R_{4n+1}=\lambda
\\
R_{4n+2}+R_{4n+3}=\lambda
\end{split}
\end{equation}
Also, from the recurrent formula \eqref{r3} it follows that:
\begin{equation}\label{mulitplicativeimplies}
\begin{split}
R_{8n}R_{8n-1}=R_{4n}
\\
R_{8n+2}R_{8n+1}=R_{4n+1}
\end{split}
\end{equation}
\begin{equation*}
R_{8n+4}R_{8n+3}=R_{4n+2}
\end{equation*}
\begin{equation*}
R_{8n+6}R_{8n+5}=R_{4n+3}
\end{equation*}
By combining \eqref{additiveimplies} with \eqref{mulitplicativeimplies} we get:
\begin{equation*}
R_{8n}(\lambda-R_{8n-2})=R_{4n}
\end{equation*}
\begin{equation*}
R_{8n+2}(\lambda-R_{8n})=\lambda-R_{4n}
\end{equation*}
\begin{equation*}
R_{8n+4}(\lambda-R_{8n+2})=R_{4n+2}
\end{equation*}
\begin{equation*}
R_{8n+6}(\lambda-R_{8n+4})=\lambda-R_{4n+2}
\end{equation*}
These can be transformed into the following:
\begin{equation}\label{recursivetransformed}
\begin{split}
R_{8n}=\frac{R_{4n}}{\lambda-R_{8n-2}}
\\
1-R_{8n+2}=1-\frac{\lambda-R_{4n}}{\lambda-R_{8n}}=
1-\frac{\lambda-R_{4n}}{\lambda-\frac{R_{4n}}{\lambda-R_{8n-2}}}=
\frac{R_{4n}(\lambda-1-R_{8n-2})}{\lambda(\lambda-R_{8n-2})-R_{4n}}
\\
R_{8n+4}=\frac{R_{4n+2}}{\lambda-R_{8n+2}}
\\
1-R_{8n+6}=\frac{R_{4n+2}(\lambda-1-R_{8n+2})}{\lambda (\lambda-R_{8n+2})-R_{4n+2}}
\end{split}
\end{equation}
\subsection{Step 2}
For convenience let's denote:
\begin{equation*}
\sigma=\frac{1}{\lambda-1}
\end{equation*}
Then we can rewrite the conclusion of theorem \eqref{t.main} in a concise form as follows:
\begin{equation}\label{inductivestatement}
\begin{split}
0<R_{4j}\leq \sigma
\\
0<1-R_{4j+2}\leq \sigma
\\
\lambda-1<R_{4j+3}\leq \lambda-1+\sigma
\\
\lambda-\sigma\leq R_{4j+1}<\lambda
\end{split}
\end{equation}
We will be proving that system of inequalities by induction over $j$.
Indeed, assume that \eqref{inductivestatement} holds for $j<2n$. We need to prove that
\begin{equation}\label{inductivestep1}
\begin{split}
0<R_{4\cdot 2n}\leq \sigma
\\
0<1-R_{4\cdot 2n+2}\leq \sigma
\\
0<R_{4\cdot 2n+4}\leq \sigma
\\
0<1-R_{4\cdot 2n+6}\leq \sigma
\end{split}
\end{equation}
and
\begin{equation}\label{inductivestep2}
\begin{split}
\lambda-\sigma\leq R_{4\cdot 2n+1}<\lambda
\\
\lambda-1<R_{4\cdot 2n+3}\leq \lambda-1+\sigma
\\
\lambda-\sigma \leq R_{4\cdot 2n+5}<\lambda
\\
\lambda-1< R_{4\cdot 2n+7}\leq \lambda-1+\sigma
\end{split}
\end{equation}
In fact, we are only concerned with proving \eqref{inductivestep1} 
as \eqref{inductivestep2} will then follow automatically from \eqref{r2}.
\\
\subsection{Step 3}
Applying the inequalities \eqref{recursivetransformed} and the 
inductive assumption \eqref{inductivestatement} we get:
\begin{equation}\label{firsthalf}
\begin{split}
R_{8n}=\frac{R_{4n}}{\lambda-R_{8n-2}}\leq \frac{\sigma}{\lambda-1}\leq \sigma
\\
R_{8n}=\frac{R_{4n}}{\lambda-R_{8n-2}}> 0
\\
1-R_{8n+2}=\frac{R_{4n}(\lambda-1-R_{8n-2})}{\lambda(\lambda-R_{8n-2})-R_{4n}}
\geq \frac{\sigma(\lambda-2)}
{\lambda (\lambda-1+\sigma)}> 0
\\
1-R_{8n+2}=\frac{R_{4n}(\lambda-1-R_{8n-2})}{\lambda(\lambda-R_{8n-2})-R_{4n}}
\leq \frac{\sigma(\lambda-2+\sigma)}{\lambda (\lambda-1)-\sigma}\leq \sigma
\end{split}
\end{equation}
The very last inequality in \eqref{firsthalf} follows from the following observation:
\begin{equation*}
\frac{\sigma(\lambda-2+\sigma)}{\lambda(\lambda-1)-\sigma}\leq\sigma
 \Leftrightarrow
\frac{\lambda-2+\sigma}{\lambda(\lambda-1)-\sigma}\leq 1
\Leftrightarrow
\lambda-2+2\sigma\leq \lambda(\lambda-1)
\Leftrightarrow
\lambda-2+\frac{2}{\lambda-1}\leq \lambda(\lambda-1)
\Leftrightarrow
\end{equation*}
\begin{equation*}
\Leftrightarrow
2\leq (\lambda-1) (\lambda^2-2\lambda+2)
\Leftrightarrow
2=min_{2\leq \lambda} (\lambda-1)((\lambda-1)^2+1)
\end{equation*}
Thus \eqref{firsthalf} proves the first two inequalities of \eqref{inductivestep1}.
As for the other two inequalities of \eqref{firsthalf}
we apply the inequalities \eqref{recursivetransformed} and the 
inductive assumption \eqref{inductivestatement} to get:
\begin{equation}\label{secondhalf}
\begin{split}
R_{8n+4}=\frac{R_{4n+2}}{\lambda-R_{8n+2}}\leq \frac{1}{\lambda-R_{8n+2}}
\\
1-R_{8n+6}=\frac{R_{4n+2}(\lambda-1-R_{8n+2})}{\lambda (\lambda-R_{8n+2})-R_{4n+2}}\leq 
\frac{\lambda-1-R_{8n+2}}{\lambda (\lambda-R_{8n+2})-1}
\end{split}
\end{equation}
By inserting the inequalities that we obtained for $R_{8n+2}$ in \eqref{firsthalf} 
into \eqref{secondhalf} we obtain:
\begin{equation}\label{secondhalf2}
\begin{split}
R_{8n+4}=\frac{1}{\lambda-R_{8n+2}}\leq \frac{1}{\lambda-1}
\\
1-R_{8n+6}=\frac{\lambda-1-R_{8n+2}}{\lambda (\lambda-R_{8n+2})-1}
\leq \frac{\lambda-2+\sigma}{\lambda(\lambda-1)-1}\leq \sigma
\end{split}
\end{equation}
The very last inequality in \eqref{secondhalf2} follows from the following observation:
\begin{equation*}
\frac{\lambda-2+\sigma}{\lambda(\lambda-1)-1}\leq \sigma
 \Leftrightarrow
 \lambda-2+\sigma\leq \sigma (\lambda(\lambda-1)-1)
 \Leftrightarrow
\end{equation*}
\begin{equation*}
 (\lambda-2)(\lambda-1)+1\leq \lambda(\lambda-1)-1
 \Leftrightarrow
 4\leq 2\lambda
\end{equation*}

\end{document}